\documentclass[11pt]{amsart}
\usepackage{amsmath,amssymb, amsthm}

\newcommand\R{{\mathbb R}}

\newcommand\Z{{\mathbb Z}}
\newcommand\Q{{\mathbb Q}}
\newcommand\semi{\propto}

\newtheorem{theorem}{Theorem}[section]
\newtheorem{proposition}[theorem]{Proposition}
\newtheorem{lemma}[theorem]{Lemma}
\newtheorem{corollary}[theorem]{Corollary}

\begin{document}

\title{The Large Scale Geometry of the Higher Baumslag-Solitar Groups}
\author{Kevin Whyte}

\maketitle

\section*{Introduction}
  The Baumslag-Solitar groups: 

$$BS(m,n)=<x,y| x y^{m} x^{-1} = y^{n}>$$

are some of the simplest interesting infinite groups which are not 
lattices in Lie groups.  They have been studied in depth from the  
point of view of combinatorial group theory.  It is natural to ask 
if the geometric approach to the theory of infinite groups, which has been so 
successful in the study of lattices, can yield any insights in this 
nonlinear case.

  The first step towards a geometric understanding of the Baumslag-Solitar 
groups is to decide which among the $BS(m,n)$ are quasi-isometric.  The 
groups $BS(1,n)$ are solvable, hence amenable, and so are not quasi-isometric 
to any of the $BS(m,n)$ with $1 < m \leq n$ which contain free subgroups
and hence are are nonamenable.

  The solvable groups $BS(1,n)$ are in many respects the 
most lattice-like of the Baumslag-Solitar groups.  They are discrete
subgroups in products of real and $p$-adic Lie groups.  
The groups $BS(1,n)$ are classified up to quasi-isometry by Farb and 
Mosher in \cite{FM1}.  They prove that $BS(1,n)$ and $BS(1,m)$ 
are quasi-isometric only if $n$ and $m$ have common powers.  When $n$ 
and $m$ have common powers $BS(1,n)$ and $BS(1,m)$ are not only 
quasi-isometric, but are commensurable (have isomorphic subgroups of 
finite index).  This is the same rigidity phenomenon as occurs
for nonuniform lattices in higher rank.  Despite this rigidity, their
full group of self quasi-isometries is quite large, and in this they
more closely resemble uniform lattices.

  In this paper we classify all the Baumslag-Solitar groups up to 
quasi-isometry.  The higher Baumslag-Solitar groups, namely those with 
$1<m<n$, are unlike the groups $BS(1,n)$ in many ways.  They are nonlinear,
not residually finite, and usually not Hopfian.  Indeed, this ``bad'' behavior
was the motivation for their discovery.  Our results show that the higher 
Baumslag-Solitar groups exhibit a surprising lack of rigidity;  all
the higher Baumslag-Solitar groups,  aside from the degenerate case 
of $BS(n,n)$, are quasi-isometric to each other.  The quasi-isometries
we construct do not reflect any clear algebraic relationship between
the groups.  In particular, many of the groups we prove to be quasi-isometric 
are not commensurable.  Our method of constructing quasi-isometries 
seems to be fundamentally different from earlier constructions.

  Our main results concern a class of groups somewhat larger than the 
class of Baumslag-Solitar groups.  We define a {\em graph of ${\mathbb Z}$s} 
as a finite graph of groups, in the sense of Serre (\cite{Se}), with 
all vertex and edge groups infinite cyclic.  This class includes the 
Baumslag-Solitar groups, which are precisely the HNN extensions of $\Z$. 

\begin{theorem}[Classification of Graphs of ${\mathbb Z}$s]\label{gen} If $G$ 
is a graph of ${\mathbb Z}$s and $\Gamma=\pi_{1} G$ then exactly one of the 
following is true:

\begin{enumerate}
\item $\Gamma$ contains a subgroup of finite index of the form 
$F_{n} \times {\mathbb Z}$.
\item $\Gamma = BS(1,n)$ for some $n>1$.
\item $\Gamma$ is quasi-isometric to $BS(2,3)$.
\end{enumerate}
\end{theorem}

Here $F_{n}$ is the free group on $n$ generators.

\begin{corollary}[Classification of Higher Baumslag-Solitar Groups]\label{bs} All 
the groups $BS(m,n)$ with $1<m<n$ are quasi-isometric to each other.
\end{corollary}

\subsection{Outline}

Let $G$ be the fundamental groups of a graph of $\Z$s.  We begin by 
constructing a geometric model space for the geometry of $G$.  This 
model is a contractible $2$-complex, $X_{G}$ on which $G$ acts cocompactly, 
freely and properly discontinuously by isometries.  The algebraic fact 
that $G$ is a graph of $\Z$s translates into the geometric fact that 
$X_{G}$ is a warped product of a tree with $\R$.  In other words, 
$X_{G}$ is topologically $T \times \R$, with a metric which differs 
from the product metric in that the metric on $v \times \R$ is scaled 
by a warping function $T \to \R^{+}$.

The tree $T$ is the Bass-Serre tree of the graph of groups, and the 
warping function is induced by a $G$ invariant orientation on $T$.  If
two graphs of $\Z$s are quasi-isometric, we show that there is a 
quasi-isometry between their Bass-Serre trees which coarse respects, 
in an appropriate sense, the orientations.  Conversely, any coarsely
orientation preserving quasi-isometry between Bass-Serre trees induces
a quasi-isometry between the groups.  Thus the classification of 
Baumslag-Solitar groups, and graphs of $\Z$s in general, reduces to 
classifying coarsely oriented trees.

The heart of our construction is the construction of coarsely 
orientation preserving quasi-isometries between trees.  We first
decompose the trees into lines of constant ``slope'' (see \S 
\ref{laminate}) with respect to the orientation.  The quasi-isometries
are built line by line.  This also requires a quasi-isometry between
the spaces of lines, with nice properties with respect to the 
orientation.  Building this matching of lines uses Hall's Selection
theorem, and the axiom of choice.   

This construction is sufficient to allow us to completely classify
graphs of $\Z$s up to quasi-isometry.  We also explore the issue
of commensurabilities among graphs of $\Z$s sufficiently to show
that although all the groups $BS(m,n)$ for $1<m<n$ are 
quasi-isometric, they are in general not commensurable.  Thus we have 
many new explicit examples of groups which are quasi-isometric but not
commensurable.

We next turn to describing the quasi-isometry groups of these graphs
of $\Z$s.  For the quasi-isometry groups of the solvable Baumslag-Solitar 
groups, $BS(1,n)$, there is a nice description in \cite{FM1}.  The 
situation for the higher Baumslag-Solitar groups is substantially 
more complex.  We discuss the complications and give several 
descriptions, none entirely satisfactory, of these quasi-isometry groups.

We discuss some generalizations.  It is natural to ask what an
arbitrary finitely generated groups quasi-isometric to a graph
of $\Z$s can be.  It follows from \cite{MSW} that any such group
is a finite graph of virtual $\Z$s.  As our classification
extends without change to that larger class of groups, we get a
complete description of the class of groups quasi-isometric to the
higher Baumslag-Solitar groups.  We also discuss some further classes
of graphs of groups to which our classification of coarsely 
oriented trees is relevant.

\section{The geometric models}\label{models}

  One of the basic principles of geometric group theory is the
Milnor-Svarc theorem, which says that if $G$ is a finitely generated
group, and $G$ acts properly discontinuously and cocompactly by 
isometries on a proper geodesic metric space $X$, then $X$ is 
quasi-isometric to $G$.  Thus, for questions about the quasi-isometric
geometry of $G$, one can work instead with $X$.

\subsection{The $2$-complexes}

  Let $\Gamma$ be a graph of ${\mathbb Z}$s, $G=\pi_{1}\Gamma$, and 
$T$ the Bass-Serre tree of $G$.  

  We first describe a 2-complex $X_{G}$ on which $G$ acts 
properly discontinuously and cocompactly by isometries.  Build a 
compact complex with $\pi_{1}=G$ out of the graph $\Gamma$ as
follows:  start with a disjoint collection of circles, one for 
each vertex of $\Gamma$.  For each edge of $\Gamma$ glue in an 
$S^{1}\times [0,1]$ where the attaching maps at each end are covering
maps inducing the same map on fundamental groups as the inclusions of 
the corresponding edge groups.  The universal cover of the complex is 
the desired $X_{G}$.

  Following \cite{FM1}, we give another description of $X_{G}$.  
Topologically, $X_{G}$ is $T \times \R$.  Let $e$ be an edge 
of $\Gamma$, and let the index of the inclusion into its vertex groups be 
$n \geq m$.  The action of the edge group of $e$ on the strip $e 
\times \R$ is translation by $n$ over one endpoint and by 
$m$ over the other.  This becomes isometric if we metrize the strip 
as a warped product $dt^{2} + (\frac{n}{m})^{2t}ds^{2}$, where $t$ is 
the parameter along $e$, and $x$ along $\R$.  This makes $e 
\times \R$ isometric to a horostrip (the region between two 
concentric horoballs) of width 1 in a space of constant curvature $- 
\ln \frac{n}{m}$.

For any vertex $v$ in $T$, we call the subspace $v \times \R$ of 
$X_{G}$ the {\em vertex space} over $v$.  Likewise, for any edge $e$, 
$e \times \R$ is the {\em edge space} over $e$.  

Given any two vertices of $T$, $v_{1}$ and $v_{2}$, let $G_1$ and 
$G_{2}$ be their stabilizers.  Let $G_{12}$ be their intersection, 
which is the stabilizer of the path between them.  $G_{12}$ has finite
intersection in both $G_{1}$ and $G_{2}$.  We call the ratio 
$\frac{[G_{12}:G_{1}]}{[G_{12}:G_{2}]}$ the {\em contraction factor}
between $v_{1}$ and $v_{2}$.  The terminology is justified by the
geometric interpretation as the contraction factor of the closest point 
projection map from the vertex space over $v_{1}$ to the vertex space
over $v_{2}$.  

It is more convenient to work with an additive rather than 
multiplicative invariant.  We define the {\em height change} between
two vertices as the logarithm of the contraction factor.  By
choosing a base point of $T$, we can define the {\em height} of a
vertex $v$ as the height change between the base point and $v$. We
extend the height function, $h$, to all of $T$ by linear interpolation
along edges.  

The metric on $X_{G}$ can be described in terms of the height function
as a warped product $T \semi \R$ with warping function $e^{-h}$.  Thus
$T$, together with the height function, determines the complex $X_{G}$
up to isometry.

\subsection{The coarsely oriented Bass-Seree tree}

We view the height change along edges as giving a quantitative 
analogue of an orientation.  If $S$ is an oriented tree, then we can
define a height change function by declaring the height change across 
an edge to be $1$ or $-1$ depending on the orientation.  With this
height change function, an isometry of $S$ preserves the
orientation if and only if it preserves the height change function.  
Just as we view a quasi-isometry as a ``large scale isometry'', we
view a quasi-isometry which preserves the height change function, on
a large scale, as coarsely orientation preserving.

{\bf Definition} A quasi-isometry $f:T_{1} \to T_{2}$ between trees 
with height functions $h_{1}$ and $h_{2}$, is {\em coarsely 
orientation preserving} is there is $C>0$ so that for all $v_{1}$ 
and $v_{2}$ in $T_{1}$:
$$|h_{1}(v_{1},v_{2}) - h_{2}(f(v_{1}),f(v_{2}))| \leq C$$

 Notice that only the height change between two points is involved in 
this definition, so the notion of coarsely orientation preserving is 
independent of the choice of base points.

\begin{theorem}\label{extend} For $i=1,2$, let $G_{i}$ be a graph 
of ${\mathbb Z}$s, with Bass-Serre tree $T_{i}$.  If $f$ is a 
quasi-isometry from $T_{1}$ to $T_{2}$ which is coarsely orientation 
preserving then $f \times Id : X_{G_{1}} \to X_{G_{2}}$ is a quasi-isometry.
\end{theorem}
\begin{proof}

The cases of bounded and unbounded height functions are fundamentally
different.  In the former, the vertex spaces are isometrically 
embedded, while in the latter they are exponentially distorted.  
Clearly, boundedness of height function is a coarse orientation 
preserving quasi-isometry invariant.  The theorem holds in this case 
as, when the height function is bounded, the complex $X$ is 
bilipschitz equivalent to the product $T \times \R$.  Thus we assume
the height functions are both unbounded.

\begin{lemma}[Metric Approximation] If the height function on $T$ is 
unbounded then the distance in $X_{G}$ is quasi-isometric to: 

  $$d_{T}(t_{1},t_{2}) + max(0, -h(t_{1},t_{2}) + \ln |x_{1}-x_{2}|)$$

where $h(t_{1},t_{2})$ is the maximum height along the geodesic 
$t_{1}t_{2}$.
\end{lemma}

\begin{proof}

As the set of vertex spaces is coarsely dense in $X$, we may assume
that $t_{1}$ and $t_{2}$ are vertices of $T$.  Given any path $p$ from
$(t_{1},x_{1})$ to $(t_{2},x_{2})$, we can replace $p$ by a path which
is piecewise horizontal (constant $\R$ coordinate) or vertical 
(constant $T$ coordinate) without multiplying the height by more that
a constant factor.  

The total length of such a path is the length of its projection to $T$ 
(= the length of the horizontal segments) plus the length of the 
lengths of the vertical segments.  The length of a vertical segment in 
the vertex space over $t$ is $e^{-h(t)}$ times the change in the $\R$ 
coordinate.  Thus any path can be shortened by moving all the vertical 
changes to occur in the vertex space over the point of maximal height on 
the projection of the path to $T$.  Thus the distance between the 
points $(t_{1},x_{1})$ and $(t_{2},x_{2})$ is bounded below by a multiple 
of the minimal length of a path which is a horizontal path from 
$t_{1}$ to a point $t$, followed by a vertical path from $(t,x_{1})$ to
$(t,x_{2})$ an then a horizontal path from $t$ to $t_{2}$.  The length 
of such a path is $d_{T}(t_{1},t) + d_{T}(t_{2},t) + e^{-h(t)}|x_{1}-x_{2}|$.

This length is equal to $d_{T}(t_{1},t_{2}) + 2d_{T}(t,t_{1}t_{2}) + 
e^{-h(t)}|x_{1}-x_{2}|$.  Replacing $t$ by the closest point at the
same height as $t$ to $t_{1}t_{2}$ shortens the path, so we may assume
that $t$ is this closet point.   

Since there is a cocompact symmetry group, it is easy to see that there are 
$\beta > 0$ and $C >0$ so that the distance of any point $t$ of $T$ to the 
set of points of height at least $h$ in $T$ is  within $C$ of $\beta 
max(0,h-h(t))$.  Thus the minimal length of a path from 
$(t_{1},x_{2})$ to $(t_{2},x_{2})$ is, to within $C$, the minimum over $h$ of:

$$d_{T}(t_{1},t_{2}) + 2max(0,\beta|h-h(t_{1},t_{2})|) + e^{-h}|x_{1}-x_{2}|$$

  The minimum of this over all $h$ occurs at $h=h(t_{1},t_{2})$ if 
$|x_{1} - x_{2}| \leq \frac{2e^{h(t_{1},t_{2})}}{\beta}$, and at
$h=\ln \frac{\beta |x_{1}-x_{2}|}{2}$ otherwise.  Substituting this
value for $h$ finishes the proof of the lemma.
    
\end{proof}

Using this lemma, we complete the proof of the theorem. By choosing basepoints 
so that $f$ is basepoint preserving, we may assume the difference 
$|h_{1}(t)-h_{2}(f(t))|$ is bounded.  Since $f$ is a quasi-isometry, 
the image of the geodesic from $t$ to $t'$ is within a uniformly bounded 
distance of the geodesic from $f(t)$ to $f(t')$.  Combining these facts, we 
see that the difference between $h_{1}(t,t')$ and $h_{2}(f(t),f(t'))$ is 
uniformly bounded.  The approximation to the distance in $X$ in the lemma is 
thus quasi-preserved by $f \times Id$, and thus $f \times Id$ is a quasi-isometry.

\end{proof}

Thus, if we can construct a coarse orientation preserving 
quasi-isometry between the Bass-Serre trees of two graphs of $\Z$s, 
this gives a quasi-isometry between the groups.  In fact, all 
quasi-isometries among graphs of $\Z$s arise this way, see \S 
\ref{qigroups}.  

\section{Constructing quasi-isometries}\label{trees}

{\bf Definition} A coarsely oriented tree is {\em homogeneous} if the 
multiset of height changes of edges incident to a vertex $v$ is the 
same for all $v$.  This is equivalent to the transitivity of height change 
preserving isometries.

The Bass-Serre trees of the Baumslag-Solitar groups are homogeneous.  
We show in section \S \ref{classify} that any coarsely oriented tree with 
cocompact symmetry group is coarsely orientation preserving 
quasi-isometric to a homogeneous tree.  

In this section we classify homogeneous coarsely oriented tree up to 
coarsely orientation preserving quasi-isometry.  Recall that is $T$ is
an oriented tree, there is an induced coarse orientation in which the 
height change across an edge is either $1$ or $-1$ depending on 
whether is edge is crossed with or against the orientation.  
Homogeneous oriented trees are determined by their {\em type}, which 
is the ordered pair $(n,m)$ of the number of edges oriented away from 
and the number of edges oriented towards any vertex.

\begin{theorem}[Classification of Homogeneous Trees] \label{homo} Let 
$T$ be a homogeneous coarsely oriented tree with height function $h$.  
Precisely one of the following holds:
\begin{itemize}
\item $h$ is constant.
\item At every vertex of $T$ there is one edge which strictly 
increases (resp. decreases) height, and all the other edges at the 
vertex strictly decrease (resp. increase) height.
\item $T$ is coarsely orientation preserving quasi-isometric to the 
oriented tree of type $(2,2)$.
\end{itemize}
\end{theorem}

\begin{proof}
    
The bulk of the proof of this theorem is constructing coarsely 
orientation preserving quasi-isometries to show

\begin{lemma}[The Main Lemma]\label{main} If $T$ is a homogeneous coarsely
oriented tree for  which, at every vertex, there are at least two edges
which strictly  increase height, and two edges which strictly decrease
height, then
$T$ is coarsely orientation preserving quasi-isometric to the 
homogeneous oriented tree of type $(2,2)$.
\end{lemma}

Assuming the lemma we complete the proof of the theorem.

If there are no edges which change height, then $h$ is constant.  
Otherwise there are, at every vertex, both an edge which strictly 
increases height and an edge which strictly decreases height.

Suppose there are edges which do not change height.  Let $F$ be
the forest of such edges.  The components of $F$ are either edges, 
with one at every vertex, or infinite trees without valence one
vertices.  In the former case, collapsing $F$ is a coarse orientation 
preserving quasi-isometry to a tree which satisfies the hypothesis
of lemma \ref{main}.  In the latter case, the following lemma 
produces a subset $F$, the collapsing of which has the same result.

\begin{lemma} Let $S$ be an infinite tree without valence one 
vertices.  There is a subset of the edges of $S$ which contains 
exactly one edge at every vertex.
\end{lemma}

\begin{proof}
Let $S'$ be a maximal subtree of $S$ for which there is such a 
subset of edges.  If $S' \neq S$ then there is a vertex $v$ of $S-S'$
and an edge $e$ with one endpoint $v$ and the other endpoint, $u$, in $S'$.
If $u$ is not in an edge of the subset, then one can extend $S'$ to
$S' \cup e$ and add $e$ to the subset of edges.  If $u$ is in one
of the subset of edges of $S'$ then let $e'$ be any edge at $v$ other 
than $e$, and extend $S'$ to $S' \cup e \cup e'$ adding $e'$ to the
subset of edges.  In either case this contradicts maximality.
\end{proof}    

Finally, if there are no edges which do not change height, then one
is clearly either in the second case of the theorem or satisfy the
hypotheses of lemma \ref{main} and hence in the third case.  This
completes the proof of the theorem, assuming lemma \ref{main}.

\end{proof}

\begin{proof}

We now turn to the proof of lemma \ref{main}.  There are two steps 
in this proof.  The first step is to decompose the tree into lines
along which the height function changes at essentially a constant rate
with respect to length.  The second step is to find a matching of
the lines in one tree with the lines in the other so that we can
assemble a coarsely orientation preserving quasi-isometry line by
line.

\subsection{Constant slope laminations}\label{laminate}

{bf Definition} Let $\beta$ and $C$ in ${\mathbb R}$ be given.  We call 
a bi-infinite geodesic $\gamma$ in $T$ a {\em line of slope $(\beta,C)$} if 
and only if for all $n$ and $m$ in ${\mathbb Z}$ 

$$|h(\gamma(n))-h(\gamma(m)) - \beta (n-m)| \leq C$$  

\begin{theorem}[Existence of constant slope laminations]\label{cover} 
If $T$ is a homogeneous tree with height function  which has at each vertex 
at least two edges along which the height increases, and two along which it 
decreases, then there is $\beta_{0} >0$ so that for any $0 \leq \beta \leq \beta_{0}$ 
there is a $C$ and a family of lines of slope $(\beta, C)$ exactly one of 
which passes through each vertex of $T$.
\end{theorem}

We will call such a collection a {\em lamination by lines of slope $\beta$}.
  
\begin{proof}

  Take $\beta_{0}$ to be such that there are two or more edges, at each 
vertex, which increase height by at least $\beta_{0}$ and two or more which 
decrease it by at least $\beta_{0}$.  Fix $0 \leq \beta \leq 
\beta_{0}$.  Let $M$ be the maximal amount height changes along any 
edge, and take $C=2M$.

  Given any vertex $v$ and edge $e$ at $v$ which increases height by 
at least $\beta_{0}$ we can find a ray of slope $(\beta,M)$ starting at $v$ 
and beginning with $e$.  We build this ray inductively.  If a ray $r$ of 
length $n$ has been constructed, extend it to length $n+1$ by choosing 
an edge which increases height by at least $\beta$ if $\beta (n+1) 
\geq h(r(n))-h(v)$ and choosing one which decreases height by at least $\beta$ 
otherwise.  It is easy to see that $r$ has the desired properties.  It 
is likewise possible to build a ray of slope $(-\beta,M)$ through any 
edge $e'$ at $v$ along which height decreases by at least $\beta$.  
By gluing the two we get a line of slope $(\beta,C)$.

  Now suppose we have $T'$ a subtree of $T$ which has been given a 
covering by lines of slope $(\beta,C)$.  If $T' \neq T$ then there is 
a $v$, a vertex of $T$, which is adjacent to $T'$.  Since only one edge connects $v$ to $T'$ we 
can build a line of slope $(\beta,C)$ through $v$ disjoint from 
$T'$. Then we can enlarge $T'$ to include $v$, the edge connecting $v$ to 
$T'$, and the new line.  Continuing in this way we cover all of $T$.
\end{proof}

\subsection{Matching the lines}

  Given two trees, $T_{1}$ and $T_{2}$, covered by lines of slope 
$\beta_{1}$ and $\beta_{2}$ we try to find an coarsely orientation preserving 
quasi-isometry from $T_{1}$ to $T_{2}$ one line at a time.  Given two 
lines there is an coarsely orientation preserving quasi-isometry between the 
lines, which is unique up to bounded distance.  Given a bijection between the 
sets of lines covering $T_{1}$ and those covering $T_{2}$ we get almost 
orientation preserving maps $T_{1} \to T_{2}$ and $T_{2} \to T_{1}$ with 
compositions at bounded distance from the identity maps of $T_{1}$ and $T_{2}$.  
We now discuss the precise conditions which make this map a 
quasi-isometry of the trees.

  Let $T_{1}'$ and $T_{2}'$ be the trees obtained from $T_{1}$ and 
$T_{2}$ by collapsing the lines of the laminations to points.   Suppose 
we have a tree isomorphism, $f$, between these quotients.  This gives, 
as above, $\hat{f}:T_{1} \to T_{2}$.  This $\hat{f}$ has bounded stretch 
along the lines of the laminations and is coarsely orientation preserving.  If 
we have an edge $e$ at height $h$ in $T_{1}$ which connects two lines, 
$a$ and $b$, it maps to an edge of $T_{1}'$ and so its image under 
$\hat{f}$ maps to an edge of $T_{2}'$.  There is a unique edge $e'$ 
of $T_{2}$ which maps to edge of $T'_{2}$.  The edge $e'$ connects two 
lines, $a'$ and $b'$, in $T_{2}$.  

Since $\hat{f}$ is coarsely orientation preserving, the end points of 
$e'$ are near the points of height $h$ on $a'$ and $b'$.  If $e'$ is 
at height $h'$ then these points are at distance $2|h'-h| + 1$ in $T_{2}$.  
Thus $\hat{f}$ is a quasi-isometry of $T_{1}$ and $T_{2}$ if, for 
every edge $e$ of $T'_{1}$, the heights of the edge in $T_{1}$ 
mapping to $e$ and the edge of $T_{2}$ mapping to $f(e)$ differ by
a uniformly bounded amount.

\begin{proposition} For $i=1,2$  let $T_{i}$ be a homogeneous tree of 
valence $n_{i}$ and be covered by lines of slopes $\beta_{i}$.  If 
$\frac{\beta_{1}}{\beta_{2}} = \frac{n_{1}-2}{n_{2}-2}$ then there are 
a $K>0$ and a tree isomorphism between $T_{1}'$ and $T_{2}'$ so that 
corresponding edges, when lifted to $T_{1}$ and $T_{2}$, differ in 
height by at most $K$.
\end{proposition}
\begin{proof}

Pick base points $v_{1}$ in $T_{1}'$ and in $T_{2}'$, let 
$f(v_{1})=v_{2}$.  Assume we can biject the edges at $v_{1}$ and $v_{2}$ 
in such a way as to change heights by at most $K$.  This gives $f$ on
the balls of radius $1$ around the basepoints.  Suppose we have the map
$f$ defined between the balls of radius $n$.  If, for each $v$ in the
sphere of radius $n$, we can biject the edges at $v$ which connect 
to the sphere of radius $n+1$ with those at $f(v)$ which connect to the
sphere or radius $n+1$, in such a way as to change height by at 
most $K$, then we can extend $f$ to the balls of radius $n+1$.  Then, 
by induction, we would have the desired $f$.

To construct the edge bijections needed in this construction we use 
Hall's selection theorem. In this context this says that bijections will 
exist between the edges at $w_{1}$ and $w_{2}$ if and only if for every 
interval $[a,b]$ in ${\mathbb R}$ the number of edges at $w_{1}$ with 
heights in $[a,b]$ is no more than the number at $w_{2}$ with heights 
in $[a-K,b+K]$, and vice versa.  This holds because the number of 
vertices on a line $l$ of slope $\beta$ with heights in the range $[a,b]$ is, 
to within a uniform additive error, $\frac{b-a}{\beta}$.  Thus, by the 
condition on the slopes, Hall's theorem applies for $K$ large enough. 
\end{proof}

Lemma \ref{cover} gives laminations of $T_{1}$ and $T_{2}$ by lines of
constant slope, with slopes of arbitrary ratio.   This completes the
proof of lemma \ref{main}.
\end{proof}

\section{The classification of graphs of $\Z$s}

Using theorem \ref{extend} and the previous construction, we can 
construct quasi-isometries between many graphs of $\Z$s.  For the
solvable Baumslag-Solitar groups, the classification in \cite{FM1}
proves that quasi-isometry implies abstract commensurability.  The
quasi-isometries we construct are very different in nature, relying 
on the axiom of choice.  We investigate when these graphs of $\Z$s 
are commensurable in sufficient detail to see that many of the groups
we prove are quasi-isometric are not commensurable.  In particular,
while all of the higher Baumslag-Solitar groups $BS(m,n)$ for $1<m<n$ are
quasi-isometric, they are, in general, not commensurable.

\subsection{The quasi-isometric classification}\label{classify}

\begin{proof}
Theorem \ref{homo} allows us to construct the quasi-isometries we need
to prove Theorem \ref{gen}.  We first show that if $G$ is any graph of
$\Z$s then its Bass-Serre tree is coarsely orientation preserving 
quasi-isometric to a homogeneous tree.

We can assume that there are no edges in the graph of groups, 
$\Gamma$, which have distinct endpoints and for which the edge group includes 
isomorphically to either of its vertex group.  If there were any such 
edges, they could be collapsed to give a graph of groups with the
same fundamental group and fewer edges.

\begin{lemma} Let $F$ be a maximal tree in $\Gamma$.  There is a family of 
lifts of $F$ to $T$ so that every vertex of $T$ is contained in exactly 
one of the lifts in the family.
\end{lemma}
\begin{proof}
    
This is done exactly as in Theorem \ref{cover}.  Since every edge group of $F$ 
includes as a subgroup of index at least two in both of its vertex 
groups, for any lift of an endpoint to $T$ there are at least two 
lifts of the edge at that vertex.  If we have lifts which 
cover a subtree $T'$ of $T$ then there is a $v$ in $T$ adjacent to $T'$.  
As each edge of $F$ has more than one lift at each lift of its 
endpoints there is a lift of $F$ through $v$ disjoint from $T'$.

\end{proof}

Pick a base point in $\Gamma$ and define the height of a lift of $F$ 
as the height of the lift of the base point it contains.  Then the 
tree $\hat{T}$ of these lifts, or equivalently the tree obtained by collapsing 
each lift, is a homogeneous tree coarsely orientation preserving 
quasi-isometric to $T$.  

If $T$ has bounded height function then it easy to see that $G$ has a 
subgroup of finite index which is $F_{n} \times {\Z}$.  

If the height function on $T$ is unbounded then the height function 
on $\hat{T}$ is also unbounded.  Each vertex then must have at least
one edge increasing height and one decreasing height.  If $F$ contains
any edges then there is at least one edge which does not change height 
at each vertex of the collapsed tree.  As in \S \ref{trees} this 
implies that $T$ is coarsely orientation preserving quasi-isometric
to the oriented tree of type $(2,2)$.

If $F$ contains no edges, then $\Gamma$ has only one vertex.  If there 
is a loop in $\Gamma$ which does not change height then again $T$ is coarsely 
orientation preserving quasi-isometric to the oriented tree of type $(2,2)$.
The same holds, by lemma \ref{main}, if there are two or more loops that do
change height, or a single loop which changes height which has more 
than one lift at both of its endpoints.  Thus the only graph of $\Z$s with 
unbounded height function not coarsely orientation preserving quasi-isometric 
to the oriented tree of type $(2,2)$ is a graph of $\Z$s with a single vertex
and a single edge which includes isomorphically at one end.  These are 
precisely the solvable Baumslag-Solitar group, which are classified up
to quasi-isometry in \cite{FM1}.

This completes the proof of Theorem \ref{gen}.
\end{proof}

\subsection{Noncommensurability}

We investigate when graphs of $\Z$s are commensurable.  While we do not
get a complete classification, we show that many of the groups we have
shown to be quasi-isometric are not commensurable.

\begin{proposition} Suppose $(a,b)=(c,d)=1$ and
$\frac{a}{b}\neq
\frac{c}{d}$, then the groups $BS(a,b)$ and $BS(c,d)$ are not
commensurable.
\end{proposition} 

Let $\Gamma$ be any graph of ${\mathbb Z}$s not quasi-isometric to 
$F_{n} \times {\mathbb Z}$ or to a solvable Baumslag-Solitar group.

An element $\gamma$ is of {\em vertex type} if and only if for any 
$\sigma \in \Gamma$ there are $n$ and $m$ nonzero for which 
$\gamma^{n}=\sigma \gamma^{m} {\sigma}^{{-1}}$.   

\begin{lemma} $\gamma$ is of vertex type if and only if $\gamma$ stabilizes a 
vertex. 
\begin{proof}

Since the tree has bounded valence any two vertex stabilizers are 
commensurable, so certainly any element which fixes a vertex is of 
vertex type.  Conversely, if $\gamma$ does not stabilize a vertex then 
it is a hyperbolic tree automorphism, so any element which conjugates one 
power of $\gamma$ to another must preserve its axis.  This can only be 
the entire group if $T$ is quasi-isometric to ${\mathbb Z}$ which does 
not happen for graphs of ${\mathbb Z}$s in this quasi-isometry class. 
\end{proof}
\end{lemma}

If $\Gamma$ is represented by a graph without any edge groups which 
include isomorphically into either vertex group, for example $BS(m,n)$ 
for $m$ and $n$ both greater than one, then no vertex stabilizer is 
contained in another.  In that case, the maximal cyclic 
subgroups of vertex type are precisely the vertex stabilizers, so the 
vertex set of $T$ is determined as a $\Gamma$ set by $\Gamma$.  The 
height function is also determined, as it is defined in terms of the 
modular homomorphism which is the ratio of indices of the 
intersections of two vertex stabilizers in each one.  It is not 
difficult to modify this to cope with loops which include 
isomorphically into one end.  There is some ambiguity in identifying 
the edges do to the possibility of sliding.

For the special case of the groups $BS(m,n)$ with $m$ and $n$ 
relatively prime and larger than one, the only finite index subgroups 
are graphs of ${\mathbb Z}$s with underlying graph a circle and all edge groups 
including as subgroups of index $m$ and $n$ in its vertex groups.  
As discussed above, we can therefore recover the number $\frac{n}{m}$ just from 
the isomorphism type of such a group.  Thus we see that $\frac{m}{n}$, 
and therefore $m$ and $n$, are commensurability invariants.  In other 
words, no two of these Baumslag-Solitar groups are commensurable.

\section{The group quasi-isometries}\label{qigroups}

In this section we calculate the quasi-isometry group of the groups 
$BS(m,n)$, for $1<m<n$.  As all these groups, and most graphs of ${\mathbb 
Z}$s, are quasi-isometric they all have the same quasi-isometry group.
The quasi-isometry groups of the solvable Baumslag-Solitar $BS(1,n)$
is the product $Bilip(\R) \times Bilip(\Q_{n})$ \cite{FM1}.  We give a
similar description of the quasi-isometry group of the higher 
Baumslag-Solitar groups, although the final form is substantially more
complicated.

We start by proving that the special form of the quasi-isometries 
we construct in \S \ref{trees} are, in fact, the general case. 
According to \cite{FM3}, if $F: X \to X$ is a quasi-isometry, there is a 
quasi-isometry $f: T \to T$ so that $\pi(F(x))=f(\pi(x))$ for $\pi$ the 
projection of $X$ to $T$.  

\begin{lemma} If $F: X \to X$ is a quasi-isometry covering $f: T 
\to T$ then $f$ is coarsely orientation preserving.
\begin{proof}

  For any $t$ in $T$, we define the fiber distance on 
$\{t\}\times {\mathbb R}$ as the induced path metric.  Since any 
quasi-isometry quasi-preserves the vertex spaces, it quasi-preserves the 
fiber distance.  In terms of the ${\mathbb R}$ coordinate this 
distance is just $e^{-h(t)}|x_{1}-x_{2}|$.  

  For any two $t$ and $t'$ in $T$, let $p:\{t\}\times{\mathbb R} \to 
\{t'\}\times {\mathbb R}$ be closest point projection.  We define the 
fiber distortion of $p$ as:

 $$ \frac{d_{F}(p(x),p(y))}{d_{F}(x,y)} $$

where $x$ and $y$ are any points on the vertex space over $t$, and 
$d_{F}$ is the fiber distance.  

  This distortion is $e^{h(t)-h(t')}$.  Closest point projection 
between the vertex spaces is preserved by the quasi-isometry, to within a 
distance determined by $d(t,t')$.  As we let $|x-y|$ go to infinity 
this additive constant has less and less effect on the  
distortion.  Thus the limit of distortion of points 
farther and farther apart is bounded above and below by multiples, 
depending only on the quasi-isometry constants of $F$, of 
$e^{h(t)-h(t')}$.    This shows that the height change $h(t)-h(t')$  
differs from $h(f(t))-h(f(t'))$ by at most some uniform additive 
error.  This is precisely the definition of coarsely orientation 
preserving.  

\end{proof}
\end{lemma}

Thus any quasi-isometry $F$ of $X$ covers an almost 
orientation preserving quasi-isometry $f$ of $T$.  According to the 
results of \S 1, $f \times Id$ is a quasi-isometry of 
$X$.  The quasi-isometry constants of $f \times Id$ may be much 
larger than those of $F$.  Even if $F$ was an isometry, $f \times Id$ 
need not be.  There is an extension which is better.  If $f$ is any 
coarsely orientation preserving map then define the height change of 
$f$, $h(f)$, as the height change between $t$ and $f(t)$ for some 
$t$ in $T$.  This change is defined up to an error determined by the $C$ in the 
definition of coarsely orientation preserving.  The map $f \times 
e^{-h(f)}$ is a quasi-isometry with constants that depend only on 
the quasi-isometry and coarsely orientation preserving constants of $f$.
 
The lemma shows that the group of quasi-isometries of $X$ splits as 
a semi-direct product of the group of coarsely orientation preserving 
quasi-isometries of $T$ and those quasi-isometries of $X$ which lie 
over the identity on $T$.  In the case of $BS(1,n)$ we can identify 
the coarsely orientation preserving quasi-isometries as $Bilip({\mathbb 
Q}_{n})$ and the quasi-isometries covering the identity as 
$Bilip({\mathbb R}$.  In this case the full quasi-isometry group is 
the product of the two.  The situation is more complicated in the case 
of $BS(m,n)$.  

A quasi-isometry, $F$, covering the identity takes each vertex space to 
itself.  For each $t$ in $T$, let $f_{t}: {\mathbb R} \to {\mathbb 
R}$ be the restriction of the quasi-isometry to the vertex space  
$\{t\} \times {\mathbb R}$.  If $F$ is an $(A,B)$ quasi-isometry 
then there are some $(A',B')$ for which $F$ restricted to each 
vertex space is an $(A',B')$ quasi-isometry with respect to fiber 
distance.  The fiber distance between $(t,x)$ and $(t,y)$ is 
$e^{-h(t)}|x-y|$, so $f_{t}$ is an $(A',B'e^{h(t)}$ quasi-isometry. 

Let $e$ be an edge of $T$ with endpoints $t$ and $t'$, where $h(t) 
\geq h(t')$.  The distance between $(t,x)$ and $(t',y)$ is 
$1+e^{-h}|x-y|$, and the distance between $(t,f_{t}(x))$ and 
$(t',f_{t'}(y))$ is $1+e^{-h}|f_{t}(x)-f_{t'}(y)|$.  Given that 
$f_{t}$ and $f_{t'}$ are quasi-isometries with constants as above, 
$F$ will be an $(A',B')$ quasi-isometry on the strip $e \times 
{\mathbb R}$ if and only if $|f_{t}(x)-f_{t'}(x)| \leq e^{h}(A'+B'-1)$.  

In summary, $F$ is a quasi-isometry covering the identity if and 
only if, for some $A,B,$ and $C$, for each $t$ in $T$, $f_{t}$ is an 
$(A,Be^{h(t)})$ quasi-isometry and, for each edge in $e$ in $T$, we 
have $d(f_{t},f_{t'})<Ce^{h}$.  

Consider the metric $d^{l}$ on $T$ where each edge has length 
$e^{h}$, where $h$ is the height of the higher endpoint.  We define 
the {\em lower boundary}, $\partial^{l}T$ as the ideal points of the metric 
completion of $T$ with respect to this metric.  The previous paragraph 
shows that a quasi-isometry covering the identity is a Lipschitz map 
from $(T,d^{l})$ to $QI({\mathbb R})$ so that the map $f_{t}$ is an 
$(A,Be^{h})$ quasi-isometry.  This gives a Lipschitz map from 
$\partial^{l}$ to $Bilip({\mathbb R})$ so that the image has 
uniformly bounded Lipschitz constants.  

\begin{lemma} Let $Bilip_{L}({\mathbb R})$ be the space of 
bilipschitz maps with bilipschitz constant at most $L$, equipped with the 
metric $d(f,g)=|f-g|_{\infty} + |f^{-1}-g^{-1}|_{\infty}$.  The  
quasi-isometries of $X$ covering the identity on $T$ are the 
bilipschitz maps $\partial^{l} T \to Bilip({\mathbb R})$ which are
contained in $Bilip_{L}$ for some $L$ for some $L$. 
\begin{proof}

We saw above that any quasi-isometry of $X$ induces such a map 
from $\partial^{l}T$ to $Bilip({\mathbb R})$.  So we need to see that 
any such map extends to a quasi-isometry of $X$, and that this 
extension is unique up to bounded distance.

For any $t$ in $T$, the distance from $t$ to $\partial^{l}T$ in the 
metric $d^{l}$ is bounded above and below by multiples of $e^{h}$.   
Let $F$ and $F'$ induce the same maps on $\partial^{l}T$.  For any $t$ 
in $T$, pick $a$ in $\partial^{l}T$ at minimal distance.  We must have 
a constant $K$ so that $d(F_{t},F_{a})\leq Ke^{h(t)}$ and the same 
for $F'_{t}$.  So $d(F_{t},F'_{t})\leq 2Ke^{h(t)}$.  As the distance 
along the vertex space over $t$ is scaled by $e^{-h(t)}$, this shows $F$ 
and $F'$ are at bounded distance.

Essentially the same argument allows us to construct an extension.  
Given a map on $\partial^{l} T$, and any $t$ in $T$, we pick any $a$ 
in $\partial^{l}T$ at minimal distance from $t$ and define $f_{t}$ to 
be equal to $f_{a}$.  For $t$ and $t'$ in $T$, and any $a$ and $a'$ 
in $\partial^{l}T$ at minimal distance from them, we know that 
$d^{l}(a,a')\leq Ke^{h(t)}+Ke^{h(t')}+d^{l}(t,t')$.  Since the map
on the lower boundary is $M$ Lipschitz 

  $$d(f_{a},f_{a'}) \leq Md^{l}(a,a') \leq MK(e^{h(t)}+e^{h(t')})+Md^{l}(t,t')$$
  
So long as $t$ and $t'$ are not equal, $d(t,t')\geq e^{max(h(t),h(t'))}$ so we 
have:
    
  $$ d(f_{a},f_{a'}) \leq M(2K+1)d^{l}(t,t')$$
  
Thus the extension is a quasi-isometry $T$ to $T$.
\end{proof}
\end{lemma} 

We can express this bilipschitz map from $\partial^{l}T$ to 
$Bilip({\mathbb R})$ differently: it can all be assembled into 
a single bilipschitz map $\partial^{l}T \times {\mathbb R}$ to itself 
which covers the identity map of $\partial^{l}T$.  

  Any coarsely orientation preserving quasi-isometry of $T$ induces a 
bilipschitz map of $\partial^{l}T$.  If $T$ has at least two edges 
decreasing height at each vertex, $\partial^{l}T$ is dense in the 
boundary of $T$, so the map on $\partial^{l}T$ determines the quasi-isometry 
up to bounded distance.  

 We have proven:
 
\begin{theorem} Let $T$ be the Bass-Serre tree of $BS(m,n)$ for 
$1<m<n$.  The group of coarsely orientation preserving of $T$ is a 
subgroup, $G$, of $Bilip(\partial^{l}T)$, and the group of quasi-isometries 
of $BS(m,n)$ is the group of bilipschitz bundle maps of 
$\partial^{l}T \times {\mathbb R}$ covering $G$.   
\end{theorem}

  It would be nice to understand which bilipschitz maps of the lower 
boundary come from coarsely orientation preserving quasi-isometries.  It 
seems likely that some sort of conformal structure should do the trick.

  There is also an upper boundary, defined as the limit points of $T$ 
with edges scaled by $e^{-h}$.  An coarsely orientation preserving 
quasi-isometry of $T$ also induces a bilipschitz map of this upper 
boundary.  As with the lower boundary, this boundary is typically 
dense in the full boundary and so a quasi-isometry is determined 
by its action on the upper boundary.  In the case of $BS(1,n)$ this 
boundary is ${\mathbb Q}_{n}$ and the bilipschitz group of the upper 
boundary is exactly the group of coarsely orientation preserving 
quasi-isometries.  

\section{Other Applications}

The results of \cite{MSW} show that any group quasi-isometric to a graph
of $\Z$s is a graph of virtual $\Z$s.  Graphs of virtual $\Z$s have 
models like those of \S \ref{models}, except that the vertex spaces 
are only quasi-isometric to $\Z$ rather than isomorphic to $\Z$.  This 
is all that we use about the vertex spaces, so our results apply in this
slightly greater generality.  

\begin{theorem} Let $\Gamma$ be a finitely generated group.  $\Gamma$ is
quasi-isometric to $BS(2,3)$ iff $\Gamma$ is a graph of virtual $\Z$s
which is neither commensurable to $F_n \times \Z$ nor virtually solvable.
\end{theorem}

More generally, the techniques of this paper can be used to study more
general graphs of groups.  One certainly needs to assume that the
Bass-Serre tress has bounded valence, which means that all of the
edge-to-vertex inclusions have finite index image.   In this case, all
of the edge and vertex groups are commensurable.   We call such a graph
of groups {\em homogeneous}.

Very little can be said in general, as one needs to understand the large
scale dynamics of isomorphisms among finite index subgroups of the vertex
groups.  One case where this is possible is graphs of groups in which
every vertex and edge groups is $\Z^n$ for some fixed $n$.  The
isomorphisms among the finite index subgroups can be represented as
elements of $SL_n(\Q)$.  

In order for the geometry to reduce to coarsely oriented trees, one needs
all of these isomorphisms to lie on a single one parameter subgroup of
$GL_n(\R)$.  The natural examples of this type are HNN extensions of
$Z^n$ along finite index subgroups.  Let $G$ be $\Z^n$, $G'$ and $G''$
finite index subgroups, and $T:G' \to G''$ an isomorphism.  Abstractly,
these HNN extensions are the groups:
$$\Gamma_T=<G,t| t^{-1} g t = Tg \hbox{, for } g \in G'>$$

We assume that at least one of the groups $G'$ or $G''$ is a proper
subgroup of $G$.  Groups of this type are studied in \cite{FM3}.  Recall
that the {\em Absolute Jordan form} of $T$ is the matrix which is the
Jordan form  of $T$ except that the values on the diagonal are the norms
of the eigenvalues rather than the eigenvalues themselves.  

\begin{theorem}\cite{FM3} Let  $\Gamma_T$ and $\Gamma_{T'}$ be as
above.  
\begin{itemize}
\item If $\Gamma_T$ and $\Gamma_{T'}$ are quasi-isometric, then for some
$\alpha \in \R^+$ the absolute Jordan forms of $T^\alpha$ and $T'$ are
equal.
\item If $G_T$ and $G_{T'}$ are solvable (which is equivalent to one of
the subgroups nonproper) then $G_T$ and $G_{T'}$ are quasi-isometric
iff there is an $\alpha \in \Q^+$ for which the absolute Jordan forms 
of $T^\alpha$ and $T'$ are equal. 
\end{itemize} 
\end{theorem}

The results of this paper allow us to complete the classification.

\begin{theorem} Let $\Gamma_T$ and $\Gamma_{T'}$ be as above.  If neither
is solvable, and the there is an $\alpha \in \R^+$ so that the absolute
Jordan forms of $T^\alpha$ and $T'$ are the same, then $\Gamma_T$ and
$\Gamma_{T'}$ are quasi-isometric. 
\end{theorem}

It is interesting that for the nonsolvable cases one has a
complete invariant of the quasi-isometry type, and a continuous family of
quasi-isometry types, while in the solvable cases one has a discrete
refinement of the invariant.  

We hope to explore more general homogeneous graphs of groups, and the
nature of their invariants, in future work.

\medskip
\noindent
Kevin Whyte\\
Dept. of Mathematics\\
University of Chicago\\
Chicago, Il\\
E-mail: kwhyte@math.uchicago.edu\\

\end{document}